\documentclass[12pt,a4paper]{amsart}
\usepackage[utf8]{inputenc}
\usepackage[T1]{fontenc}
\usepackage[in,headings]{fullpage}
\usepackage{amsbsy, amscd, amsfonts, amsmath, amsrefs, amssymb, amsthm}
\usepackage{graphicx}
\usepackage{float}
\usepackage{color}   
\usepackage[colorlinks=true,linkcolor=blue,citecolor=blue,urlcolor=blue]{hyperref}

\newtheorem{theorem}{Theorem}

\newtheorem{corollary}[theorem]{Corollary}
\newtheorem{proposition}[theorem]{Proposition}

\theoremstyle{definition}
\newtheorem{definition}[theorem]{Definition}

\theoremstyle{remark}

\newcommand{\C}{\mathbf{C}}
\newcommand{\Z}{\mathbf{Z}}

\newcommand{\R}{\mathbf{R}}
\newcommand{\N}{\mathbf{N}}
\renewcommand{\H}{\mathcal{H}}
\renewcommand{\Re}{\mathop{\mathrm{Re}}\nolimits}
\renewcommand{\Im}{\mathop{\mathrm{Im}}\nolimits}
\newcommand{\Rzeta}{\mathop{\mathcal R }\nolimits}
\newcommand{\Lzeta}{\mathop{\mathcal L }\nolimits}

\begin{document}

\title{Riemann's auxiliary function. Basic Results}
\author{J. Arias de Reyna}
\address{%
Universidad de Sevilla \\ 
Facultad de Matem\'aticas \\ 
c/Tarfia, sn \\ 
41012-Sevilla \\ 
Spain.} 




\email{arias@us.es, ariasdereyna1947@gmail.com}


\begin{abstract}
We give the definition, main properties and integral expressions of the auxiliary function of Riemann $\Rzeta(s)$. For example we prove \[\pi^{-s/2}\Gamma(s/2)\Rzeta(s)=-\frac{e^{-\pi i s/4}}{ s}
\int_{-1}^{-1+i\infty} \tau^{s/2}\vartheta_3'(\tau)\,d\tau.
\] Many of these results are known, but they serve as a reference. We give the values of $\Rzeta(s)$ at integers except at odd natural numbers. We have 
\[\zeta(\tfrac12+it)=e^{-i\vartheta(t)}Z(t),\quad \Rzeta(\tfrac12+it)=\tfrac12e^{-i\vartheta(t)}(Z(t)+iY(t)),\]
with $\vartheta(t)$, $Z(t)$ and $Y(t)$ real functions. 
\end{abstract}

\maketitle

\section{Introduction}

Siegel \cite{Siegel}*{eq.~(57)} introduces the integral 
\begin{equation}
\Rzeta(s)=\int_{0\swarrow1}\frac{x^{-s} e^{\pi i x^2}}{e^{\pi i x}-
e^{-\pi i x}}\,dx.
\end{equation}
and proves that this function is related to the  zeta function of Riemann by
\begin{equation}\label{zetaRzeta}
\zeta(s)=\Rzeta(s)+\chi(s)\overline{\Rzeta}(1-s);\qquad
\chi(s)=\pi^{s-1/2}
\frac{\Gamma\left(\frac{1-s}{2}\right)}{\Gamma\left(\frac{s}{2}\right)},\quad
\overline{\Rzeta}(s)=\overline{\Rzeta(\overline{s})}.
\end{equation}
Siegel says that essentially this relation is found in Riemann's nachlass. 

Our purpose here is to obtain some general properties of the function $\Rzeta(s)$, especially some integral and series representations. 

In Sections \ref{S:2} and \ref{S:3} we recall the results contained in \cite{Siegel} about the definition of $\Rzeta(s)$ and its connection with Riemann's zeta function. 
Section \ref{S:4} contains some integral representations and gives the values of $\Rzeta(s)$ for the integer values of $s$. In particular, we prove that $\Rzeta(-2n)=0$ for $n\in\N$ a natural number ($\N=\{1,2,3,\dots\}$).  Section \ref{S:5} gives a series representation for $\Rzeta(s)$ that is due to R. Kuzmin. Section \ref{S:6} gives integral representations in terms of theta functions. These can be found in Riemann's Nachlass but are not found in Siegel.  We have $\zeta(s)=\Rzeta(s)+\chi(s)\overline{\Rzeta}(1-s)$,  in Section \ref{S:7} we determine the class of functions that share this property with $\Rzeta(s)$. In Section \ref{S:8} we introduce an interesting real-valued real analytic function  $Y(t)$ such that
\[\Rzeta(\tfrac12+it)=\tfrac12e^{-i\vartheta(t)}(Z(t)+iY(t)).\]

\section{Riemann's integral}\label{S:2}
According to Siegel \cite{Siegel}*{eq.~(7)} Riemann considered the integral
\begin{equation}
\Phi(z,\tau)=\int_{0\nwarrow1}\frac{e^{-\pi i\tau u^2+2\pi i z u}}{e^{\pi i u}-
e^{-\pi i u}}\,du,
\end{equation}
where $0\nwarrow1$ denotes a line in the direction $e^{3\pi i/4}$ and crossing the real line at a point between $0$ and $1$. 

This integral has been considered by many authors, in particular Siegel \cite{Siegel}, Mordell \cite{Mordell} and Zwegers \cite{Zwegers}, where  easy proofs of the results I cite here can be found.

For $\Im\tau>0$ the integral $\Phi(z,\tau)$ is an entire function of $z$. For a fixed value of $z\in\C$, it is holomorphic function of $\tau$ in the upper half plane $\H=\{\tau\colon\Im \tau>0\}$.

\begin{theorem}
The function $\Phi(z,\tau)$ satisfies the two functional equations
\begin{equation}\label{E:Rec1}
\Phi(z+1,\tau)-\Phi(z,\tau)=\frac{i}{\sqrt{-i\tau}}e^{\frac{\pi
i}{\tau}\left(z+\frac{1}{2}\right)^2}.
\end{equation}
\begin{equation}\label{E:Rec2}
\Phi(z,\tau)=-e^{-2\pi iz-\pi i\tau}\Phi(z+\tau,\tau)+1.
\end{equation}
\end{theorem}

For a fixed $z$ the function $\Phi(z,\tau)$ extends analytically through the real axis and its value at rational points can be obtained by the next theorem
\begin{theorem}
For $a\in\N$ and $b\in\N$ we have
\begin{multline*}\label{E:IntegRiemann}
\bigl(1-(-1)^{b+ab}e^{-2\pi i bz}\bigr)\Phi(z,\tfrac{a}{b})=\\
\sum_{n=0}^{b-1}(-1)^ne^{-2\pi i nz}e^{-\pi i \frac{a}{b} n^2}+
(-1)^{b+ab}e^{-2\pi ibz}\frac{i}{\sqrt{-i\frac{a}{b}}}\sum_{m=0}^{a-1}e^{\pi i\frac{b}{a}(z+m+1/2)^2},
\end{multline*}
When the integer $a$ is negative $a\le -1$ the corresponding value is 
\[\sum_{n=0}^{b-1}(-1)^ne^{-2\pi i nz}e^{-\pi i \frac{a}{b} n^2}-
(-1)^{b+ab}e^{-2\pi ibz}\frac{i}{\sqrt{-i\frac{a}{b}}}\sum_{m=a}^{-1}e^{\pi i\frac{b}{a}(z+m+1/2)^2}.\]
In both cases, taking $\sqrt{-ia/b}$ with a positive real part.
\end{theorem}

\begin{corollary}
\begin{equation}\label{E:riemcomputed}
\int_{0\swarrow1}\frac{e^{\pi i u^2+ z u}}{ e^{\pi i u}-e^{-\pi i u}}\,du
=-\frac{1}{1-e^{-z}}+\frac{e^{i z^2/4\pi}}{ e^{z/2}-e^{-z/2}}.
\end{equation}
\end{corollary}

\section{Riemann auxiliary function}\label{S:3}
\begin{definition}
For $s\in\C$ we define
\begin{equation}\label{E:defRzeta}
\Rzeta(s)=\int_{0\swarrow1}\frac{x^{-s} e^{\pi i x^2}}{e^{\pi i x}- e^{-\pi i x}}\,dx,
\end{equation}
where $0\swarrow1$ denotes the line of direction $e^{-3\pi i/4}$ passing through the point $x=\frac12$, and $x^{-s}= \exp(-s\log x)$ with the main  branch of $\log x$.
\end{definition}
The function $\Rzeta(s)$ is entire. It is related to Riemann's zeta function
\begin{theorem}[Riemann-Siegel integral formula]
For any $s\in\C$, except for poles, we have
\begin{equation}\label{E:RiemSiegel}
\pi^{-\frac{s}{2}}\Gamma(\tfrac{s}{2})\zeta(s)=
\pi^{-\frac{s}{2}}\Gamma(\tfrac{s}{2})\int_{0\swarrow1}\frac{x^{-s} e^{\pi i x^2}}{e^{\pi i x}- e^{-\pi i x}}\,dx+
\pi^{-\frac{1-s}{2}}\Gamma(\tfrac{1-s}{2})\int_{0\searrow1}\frac{x^{s-1} e^{-\pi i x^2}}{e^{\pi i x}- e^{-\pi i x}}\,dx.
\end{equation}
\end{theorem}
We introduce the conjugate function $\overline{\Rzeta}(s)=\overline{\Rzeta(\overline{s})}$, and notice that
\[\overline{\Rzeta}(1-s)=\overline{\Rzeta(1-\overline{s})}
=\overline{\int_{0\swarrow1}\frac{x^{\overline{s}-1} e^{\pi i
x^2}}{e^{\pi i x}- e^{-\pi i x}}\,dx}=
\int_{0\searrow1}\frac{x^{s-1}e^{-\pi i x^2}}{ e^{\pi i x}-e^{-\pi
i x}}\,dx.\]
Hence, we may write \eqref{E:RiemSiegel} as 
\begin{equation}\label{E:Siegel}
\pi^{-\frac{s}{2}}\Gamma(\tfrac{s}{2})\zeta(s)=\pi^{-\frac{s}{2}}\Gamma(\tfrac{s}{2})\Rzeta(s)+\pi^{-\frac{1-s}{2}}\Gamma(\tfrac{1-s}{2})\overline{\Rzeta}(1-s),
\end{equation}
or equivalently 
\begin{equation}\label{E:zetaRzeta}
\zeta(s)=\Rzeta(s)+\chi(s)\overline{\Rzeta}(1-s),\qquad \chi(s)=\pi^{s-\frac12}\frac{
\Gamma(\frac{1-s}{2})}{\Gamma(\frac{s}{2})}.
\end{equation}
In the critical line $s=\frac12+it$, the two terms on the right hand side of \eqref{E:Siegel}
are complex conjugate so that we obtain
\begin{equation}\label{E:RSincrit}
\pi^{-\frac14-i\frac{t}{2}}\Gamma(\tfrac14+i\tfrac{t}{2})\zeta(\tfrac12+it)=
2\Re\{\pi^{-\frac14-i\frac{t}{2}}\Gamma(\tfrac14+i\tfrac{t}{2})\Rzeta(\tfrac12+it)\}.
\end{equation}
The phase (see \cite{AL} for the definition) of $\zeta(\frac12+it)$ is 
\begin{equation}
\vartheta(t)=\Im\log\Gamma(\tfrac14+i\tfrac{t}{2})-\tfrac{t}{2}\log\pi.
\end{equation}
Then $\vartheta(t)$ is an odd real analytic function and there exists another real analytic even function $Z(t)$, the \emph{Riemann-Siegel function}\footnote{Sometimes called Hardy's function. But it is found in Riemann's Nachlass.} such that 
\begin{equation}
\zeta(\tfrac12+it)=e^{-i\vartheta(t)}Z(t),\qquad t\in\R.
\end{equation}
And \eqref{E:RSincrit} implies 
\begin{equation}
Z(t)=2\Re\{e^{i\vartheta(t)}\Rzeta(\tfrac12+it)\},\qquad t\in\R.
\end{equation}

\section{Two integral representations and values at integers}\label{S:4}

\begin{proposition}\label{T:Rint2}
For $\Re s<0$ we have
\begin{equation}\label{E:anterior}
\Rzeta(s)=\omega e^{\pi i s/4}\sin\frac{\pi
s}{2}\int_0^{+\infty}\frac{y^{-s}e^{-\pi y^2}}{\sin\pi\omega y}\,dy,
\end{equation}
where $\omega=e^{\pi i/4}$
\end{proposition}

\begin{proof}
Start with the definition
\[\Rzeta(s)=\int_{0\swarrow1}\frac{x^{-s}e^{\pi i x^2}}{e^{\pi i
x}-e^{-\pi i x}}\,dx.\] When  $\sigma=\Re s<0$ we may substitute the line of integration 
to the parallel line through $0$. 
The opposite of this path of integration is parametrized by 
$x=\omega y$ where $\omega=e^{\pi i/4}$. For  $y>0$ we have
\begin{align*}
(y\omega)^{-s}&=\exp\{-s(\log y+\pi i/4)\}=e^{-\pi i s/4} y^{-s};\\
(-y\omega)^{-s}&=\exp\{-s(\log y-3\pi i/4)\}=e^{3\pi i s/4}
y^{-s}.
\end{align*} 
It follows that 
\[\Rzeta(s)=-\omega(e^{-\pi i s/4}-e^{3\pi i
s/4})\int_0^{+\infty}\frac{y^{-s}e^{-\pi y^2}}{2i\sin\pi\omega
y}\,dy\] Simplifying, we get  \eqref{E:anterior}.
\end{proof}

\begin{proposition}\label{defalternativa} 
For any $s\in\C$, with $s\ne2n+1$ we have
\begin{equation}\label{E:Rintegral3}
\Rzeta(s)=-\omega\frac{e^{7\pi i s/4}}{(1+e^{\pi i s})}\int_C
\frac{y^{-s} e^{-\pi  y^2}}{2i\sin\pi\omega y}\,dy,
\end{equation}
where $C$ denotes Hankel contour,  a path that start at $+\infty$ follow the real axis to $0$, encircles the origin counter clockwise and turns to $+\infty$ again, and where \[y^{-s}=\exp\{-s\log|y|-is\arg y\}\qquad\text{ with }0<\arg y<2\pi.\]
\end{proposition}

\begin{proof}
It is easy to see that when $\sigma<0$ we have
\[\int_C\frac{y^{-s} e^{-\pi
y^2}}{2i\sin\pi\omega y}\,dy
= -\int_0^{+\infty}\frac{y^{-s}
e^{-\pi  y^2}}{2i\sin\pi\omega y}\,dy+e^{-2\pi i
s}\int_0^{+\infty}\frac{y^{-s} e^{-\pi y^2}}{2i\sin\pi\omega
y}\,dy.
\]

Applying Proposition \ref{T:Rint2} this can be written as 
\[\Rzeta(s)=\omega e^{-\pi i s/4}\frac{e^{\pi i s}-1}{e^{-2\pi i
s}-1} \int_C\frac{y^{-s} e^{-\pi  y^2}}{2i\sin\pi\omega y}\,dy.\]
That simplifies to \eqref{E:Rintegral3}. 
\end{proof}

\begin{corollary}
$\Rzeta(-2n)=0$ for $n=1$, $2$, \dots\ As in the case of the zeta function, we call these trivial zeros of $\Rzeta(s)$. 
\end{corollary}
\begin{proof}
The integral in \eqref{E:Rintegral3} vanishes for $s=-2n$. It reduces to the integral on a circle and Cauchy's Theorem proves that it is equal to $0$.  For negative odd integers, the integral also vanishes, but the factor of the integral has a pole there. 
\end{proof}

\begin{corollary}
The value of $\Rzeta(s)$ at  the even integers can be given explicitly
\begin{equation}\label{E:Rzetaevenvalues}
\Rzeta(2n)=-\sum_{k=0}^n2(2^{2k-1}-1)\frac{(\pi
i)^{n-k}}{(n-k)!}\frac{\zeta(2k)}{2^{2k}}.\qquad (n=0, 1, 2,
\dots)
\end{equation}
In particular,
\begin{gather}
\Rzeta(0)=-\frac{1}{2};\quad \Rzeta(2)=-\frac{\pi^2}{12}-\frac{\pi
i}{2};\quad
\Rzeta(4)=-\frac{7\pi^4}{720}-\frac{i\pi^3}{12}+\frac{\pi^2}{4};\\
\Rzeta(6)=-\frac{31\pi^6}{30240}-\frac{7i\pi^5}{720}+\frac{\pi^4}{24}+
\frac{i\pi^3}{12}.
\end{gather}
\end{corollary}

\begin{proof}
The Taylor series of $x/\sin x$ valid for $|x|<\pi$ can be written
\begin{equation}\label{E:sincexpansion}
\frac{x}{\sin x}=1+\sum_{n=1}^\infty4(2^{2n-1}-1)\zeta(2n)
\Bigl(\frac{x}{2\pi}\Bigr)^{2n}=\sum_{n=0}^\infty4(2^{2n-1}-1)\zeta(2n)
\Bigl(\frac{x}{2\pi}\Bigr)^{2n}.
\end{equation}
The integrand in the integral expression  \eqref{E:Rintegral3}  of $\Rzeta(s)$ is given for 
$|y|<1$ by
\[\begin{aligned}
\frac{y^{-s}e^{-\pi y^2}}{2i\sin(\pi \omega y)}&=\frac{y^{-s}}{2\pi i\omega y}
\sum_{j=0}^\infty\frac{(-\pi )^j}{ j!}y^{2j}\sum_{k=0}^\infty4(2^{2k-1}-1)\zeta(2k)
\Bigl(\frac{\omega y}{2}\Bigr)^{2k}\\
&=\frac{y^{-s}}{2\pi i\omega }\frac{1}{ y}\sum_{n=0}^\infty y^{2n}
\sum_{k=0}^n4(2^{2k-1}-1)\frac{(-\pi )^{n-k}}{ (n-k)!}\frac{i^k\zeta(2k)}{2^{2k}}.
\end{aligned}\]
When $s$ is a nonnegative integer, the integral reduces to an integral along a circle of
center $0$ and arbitrary radius $<\pi$. In this case, the integrand is a uniform meromorphic function. The integral is easily computed because it is equal to the coefficient of $y^{-1}$.
We get  \eqref{E:Rzetaevenvalues}.

Notice that when $s=2n+1$, the integral is also easily computed with value $0$, but we cannot deduce the value of $\Rzeta(2n+1)$ because the factor in front of the integral has a pole there.
\end{proof}

\begin{proposition} The value of $\Rzeta(s)$ for odd negative integers can be explicitly computed.
For all $n\in\N$ 
\begin{equation}\label{E:Rzetaoddvalues}
\frac{\Rzeta(1-2n)}{(2n-1)!}\\= (-1)^n\frac{4(2^{2n}-1)\zeta(2n)}{(4\pi)^{2n}}+
4\sum_{k=0}^{n-1}(-1)^k
\frac{(2^{2k-1}-1)\zeta(2k)i^{n-k}}{(n-k)! (4\pi)^{n+k}}.
\end{equation}
In particular,
\begin{gather}
\Rzeta(0)=-\frac{1}{2};\ \Rzeta(-1)=
\Bigl(\frac{i}{\pi}-\frac{1}{2}\Bigr)\frac{1!}{2^2};\ \Rzeta(-3)=
\Bigl(-\frac{1}{\pi^2}-\frac{i}{3\pi}+\frac{1}{12}\Bigr)\frac{3!}{2^5};\\
\Rzeta(-5)=\Bigl(-\frac{i}{\pi^3}+\frac{1}{2\pi^2}+\frac{7i}{60\pi}-
\frac{1}{40}\Bigr) \frac{5!}{3\cdot2^7}.
\end{gather}
Recall that we have  $\Rzeta(-2n)=0$
for  $n\in \N$.
\end{proposition}

\begin{proof}
We get 
\begin{equation}
\int_{0\swarrow1}\frac{e^{\pi i x^2+ z x}}{ e^{\pi i x}-e^{-\pi i x}}\,dx=
\sum_{n=0}^\infty\Rzeta(-n)\frac{z^n}{n!},
\end{equation}
expanding the exponential $e^{zx}$ in the integrand, integrating term by term and applying
the definition of $\Rzeta(s)$ \eqref{E:defRzeta} of 

Recall also the integral computed by Riemann \eqref{E:riemcomputed}
\[\int_{0\swarrow1}\frac{e^{\pi i x^2+ z x}}{ e^{\pi i x}-e^{-\pi i x}}\,dx
=-\frac{1}{1-e^{-z}}+\frac{e^{i z^2/4\pi}}{ e^{z/2}-e^{-z/2}}.\]
We expand in the Taylor series. First, recall the common expansion
\begin{equation}
\frac{x}{e^x-1}=1-\frac{x}{2}-{2}\sum_{n=1}^\infty(-1)^n\zeta(2n)
\Bigl(\frac{x}{2\pi}\Bigr)^{2n}.
\end{equation}
Therefore,
\begin{equation}\label{E:exp1}
-\frac{1}{1-e^{-z}}=-\frac{1}{z}\frac{-z}{ e^{-z}-1}=
-\frac{1}{z}\Bigl(1+\frac{z}{2}-{2}\sum_{n=1}^\infty
(-1)^n\zeta(2n)\Bigl(\frac{z}{2\pi}\Bigr)^{2n}\Bigr).
\end{equation}
In the sinc expansion \eqref{E:sincexpansion} substitute $ix$ instead of $x$
\[\frac{x}{\sinh x}=4\sum_{n=0}^\infty (-1)^n (2^{2n-1}-1)\zeta(2n)\Bigl(\frac{x}{2\pi}\Bigr)^{2n}.\]
Therefore,
\begin{align}\label{E:exp2}
\frac{e^{i z^2/4\pi}}{2\sinh(z/2)}&=\frac{1}{z}
\sum_{j=0}^\infty \Bigl(\frac{i}{4\pi}\Bigr)^j\frac{z^{2j}}{j!}\sum_{k=0}^\infty4(-1)^k(2^{2k-1}-1)
\frac{\zeta(2k)}{(4\pi)^{2k}}z^{2k}\\
&=\frac{1}{z}\sum_{n=0}^\infty z^{2n}\sum_{k=0}^n
4(-1)^k(2^{2k-1}-1)\frac{\zeta(2k)}{(4\pi)^{2k}}\Bigl(\frac{i}{4\pi}\Bigr)^{n-k}\frac{1}{(n-k)!}.
\end{align}
From \eqref{E:exp1}, \eqref{E:exp2}, equating the coefficients of the two expansions of the
same function, we get \eqref{E:Rzetaoddvalues}.
\end{proof}

\section{Kuzmin representation}\label{S:5}

Kuzmin, independently of Riemann, obtained the representation \eqref{E:RiemSiegel} of the zeta function. His derivation yields another representation of $\Rzeta(s)$ as a series of incomplete Gamma functions.  This is obtained by expanding $1/\sin x$ into simple fractions and integrating term by term in the definition of $\Rzeta(s)$. 

\begin{proposition}\label{P:KuzminInt}
For  $s\in \C$ and $n\in\N$ we have
\begin{equation}\label{E:Kuzmin}
\int_{0\swarrow1}\frac{x^{2-s}e^{\pi i
x^2}}{x^2-n^2}\frac{dx}{x}=(-1)^n\frac{\pi
i}{n^s}\frac{\Gamma(s/2,\pi i n^2)}{\Gamma(s/2)}.
\end{equation}
\end{proposition}

\begin{proof}
Assume first  $\sigma=\Re s<2$, in this case we may move the line of integration 
to the parallel through the origin. If we assume also $\sigma>0$ we may apply Cauchy's
Theorem and rotate the line of integration to coincide with the imaginary axis. In fact, the
integral on the piece of circle with center at the origin and radius $R$ that joins the two paths  is bounded  by $C e^{\pi |t|}R^{2-\sigma} R^{-2} $, which tends to $0$ for $R\to+\infty$.

Then, denoting by $I$ the integral in our Proposition, we have
\begin{align*}
I&=-\int_0^{+i\infty}\frac{x^{2-s}e^{\pi i
x^2}}{x^2-n^2}\frac{dx}{x}+\int_0^{-i\infty}\frac{x^{2-s}e^{\pi i
x^2}}{x^2-n^2}\frac{dx}{x}\\
&=-e^{-\pi i s/2}\int_0^{+\infty}\frac{y^{2-s}e^{-\pi i
y^2}}{y^2+n^2}\frac{dy}{y}+e^{\pi i
s/2}\int_0^{+\infty}\frac{y^{2-s}e^{-\pi i
y^2}}{y^2+n^2}\frac{dy}{y}\\
&=2i\sin\frac{\pi s}{2}\int_0^{+\infty}\frac{y^{2-s}e^{-\pi i
y^2}}{y^2+n^2}\frac{dy}{y} \\
&=(-1)^n2i\sin\frac{\pi s}{2}\int_0^{+\infty}\frac{y^{2-s}e^{-\pi
i
(y^2+n^2)}}{y^2+n^2}\frac{dy}{y}\\
&= (-1)^n2i\sin\frac{\pi s}{2}\int_0^{+\infty}y^{2-s}\Bigl(
\int_{\pi i}^{\pi i+\infty}e^{-x(y^2+n^2)}\,dx\Bigr)\frac{dy}{y}.
\end{align*}
Now applying Fubini's Theorem,
\[I=(-1)^n2i\sin\frac{\pi s}{2}\int_{\pi i}^{\pi
i+\infty}e^{-n^2
x}\Bigl(\int_0^{+\infty}y^{2-s}e^{-xy^2}\frac{dy}{y}\Bigr)\,dx.\]
Computing the inner integral yields
\begin{align*}
I&=(-1)^ni\sin\frac{\pi s}{2}\Gamma(1-s/2)\int_{\pi i}^{\pi
i+\infty}e^{-n^2 x} x^{s/2-1}\,dx\\
&=\frac{(-1)^n\pi i}{\Gamma(s/2)}\frac{\Gamma(s/2,\pi i
n^2)}{n^s}.
\end{align*}
Both members are integer functions, so  the 
equality is true for all $s$. 
\end{proof}

\begin{theorem} For any complex number $s$
\begin{equation}\label{E:Rzetaserie}
\Rzeta(s)=-\frac{\pi^{s/2}}{s\Gamma(s/2)}e^{\pi i
s/4}+\sum_{n=1}^\infty \frac{\Gamma(s/2,\pi i
n^2)}{\Gamma(s/2)}\frac{1}{n^s},
\end{equation}
where the series is absolutely convergent.
\end{theorem}

\begin{proof}
In the definition of $\Rzeta(s)$
\[\Rzeta(s)=\int_{0\swarrow1}\frac{x^{-s}e^{\pi i x^2}}{ e^{\pi i x}-e^{-\pi i
x}}\,dx,\] we substitute the expansion of Mittag Leffler
\begin{equation}
\frac{1}{e^{\pi i x}-e^{-\pi i x}}=\frac{1}{2\pi
i}\sum_{-\infty}^{\infty}\frac{(-1)^n x}{ x^2-n^2}.
\end{equation}

For $x$ in the line of integration, this series converges absolutely to a bounded function.
We may integrate term by term with respect to the measure $x^{-s}e^{\pi i x^2}\,dx$ and the result will be an absolutely convergent series.
\[\Rzeta(s)=\sum_{n=-\infty}^{+\infty}\frac{(-1)^n}{2\pi
i}\int_{0\swarrow1}\frac{x^{2-s}e^{\pi i
x^2}}{x^2-n^2}\frac{dx}{x}.\] For $n=0$, the integral has a known value; for other values of $n$ we apply Proposition \ref{P:KuzminInt}. This yields \eqref{E:Rzetaserie}.
\end{proof}

\section{Integral representations with theta functions}\label{S:6}

We will use the classical  definitions of theta functions \cite{WW}, we consider them functions of $\tau$ taking $q=e^{\pi i\tau}$.
\begin{equation}
\begin{aligned}
\vartheta_2(\tau)&=2q^{\frac14}(1+q^2+q^6+\cdots)= \sum_{n\in\Z}q^{(n-\frac12)^2}=
2q^{\frac14}\prod_{n=1}^\infty\{(1-q^{2n})(1+q^{2n})^2\}\\
\vartheta_3(\tau)&=1+2q+2q^4+2q^9+\cdots=\sum_{n\in\Z}q^{n^2}=
\prod_{n=1}^\infty\{(1-q^{2n})(1+q^{2n-1})^2\} \\
\vartheta_4(\tau)&=1-2q+2q^4-2q^9+\cdots=\sum_{n\in\Z}(-1)^nq^{n^2}=
\prod_{n=1}^\infty\{(1-q^{2n})(1-q^{2n-1})^2\}
\end{aligned}
\end{equation}
They satisfy the transformations laws:
\begin{align}
\vartheta_2(\tau+1)&=e^{\pi i/4}\vartheta_2(\tau)& \vartheta_2(-1/\tau)&=\sqrt{-i\tau}\;\vartheta_4(\tau)\\
\vartheta_3(\tau+1)&=\vartheta_4(\tau) & \vartheta_3(-1/\tau)&=\sqrt{-i\tau}\;\vartheta_3(\tau)\\
\vartheta_4(\tau+1)&=\vartheta_3(\tau) & \vartheta_4(-1/\tau)&=\sqrt{-i\tau}\;\vartheta_2(\tau)
\end{align}
where $\sqrt{-i\tau}$ denotes the square root with a positive real part.

\begin{theorem}
Except for the poles at $s=0$ and $s=1$ 
\begin{equation}\label{E:Rthetados} 
\pi^{-s/2}\Gamma(s/2)\Rzeta(s)=-\frac{1}{s}e^{\pi i s/4}+\frac{1}{2}
e^{-\pi is/4}\int_{-1}^{-1+i\infty}\tau^{s/2}\bigl(\vartheta_3(\tau)-1\bigr)\,
\frac{d\tau}{\tau}.
\end{equation}
\end{theorem}

\begin{proof}
For $s\ne0$ or $1$, by \eqref{E:Rzetaserie} we have
\[\pi^{-s/2}\Gamma(s/2)\Rzeta(s)=-\frac{e^{\pi i
s/4}}{s}+\pi^{-s/2}\sum_{n=1}^\infty \frac{\Gamma(s/2,\pi i
n^2)}{n^s}.\] Changing variable  $x=-\pi i n^2 \tau$,
in the integral of the incomplete gamma function
\[\Gamma(s/2,\pi i n^2)=\int_{\pi i
n^2}^{\pi i n^2+\infty}
x^{s/2}e^{-x}\frac{dx}{x}=\int_{-1}^{-1+i\infty}\pi^{s/2}
n^s\Bigl(\frac{\tau}{i}\Bigr)^{s/2}e^{\pi i
n^2\tau}\frac{d\tau}{\tau}.\] Substituting into the above expression we get the following
\[\pi^{-s/2}\Gamma(s/2)\Rzeta(s)=-\frac{e^{\pi i
s/4}}{s}+\sum_{n=1}^\infty \int_{-1}^{-1+i\infty}
\Bigl(\frac{\tau}{i}\Bigr)^{s/2}e^{\pi i
n^2\tau}\frac{d\tau}{\tau}.\] Si $\tau=-1+iy$ tenemos
\[\sum_{n=1}^\infty |e^{\pi in^2\tau}|=\sum_{n=1}^\infty e^{-\pi
n^2 y}\] which for $y\to0^+$ behaves as  $\sqrt{y}$ and when
$y\to+\infty$ as $e^{-\pi y}$. This justifies the interchange of sum and integral 
\[\pi^{-s/2}\Gamma(s/2)\Rzeta(s)=-\frac{e^{\pi i
s/4}}{s}+\frac{1}{2}\int_{-1}^{-1+i\infty}
\Bigl(\frac{\tau}{i}\Bigr)^{s/2}\bigl\{\vartheta_3(\tau)-1\bigr\}
\frac{d\tau}{\tau}.\]
\end{proof}

\begin{corollary}
When $\sigma=\Re s<0$ 
\begin{equation}\label{E:Rthetauno} 
\Rzeta(s)=\frac{1}{2} e^{-\pi i
s/4}\frac{\pi^{s/2}}{\Gamma(s/2)} \int_{-1}^{-1+i\infty}
\tau^{s/2}\vartheta_3(\tau)\,\frac{d\tau}{\tau},\qquad \Re s<0. 
\end{equation}
\end{corollary}

\begin{proof}
When $\sigma<0$ the integral of 
$(\tau/i)^{s/2}\, d\tau/\tau$ can be computed explicitly. And we obtain \eqref{E:Rthetauno}.
\end{proof}

\begin{proposition}
For all $s\in\C$ we have
\begin{equation}\label{E:Rzetatres}
\pi^{-s/2}\Gamma(s/2)\Rzeta(s)=-\frac{e^{-\pi i s/4}}{ s}
\int_{-1}^{-1+i\infty} \tau^{s/2}\vartheta_3'(\tau)\,d\tau.
\end{equation}
\end{proposition}

\begin{proof}
Starting from  \eqref{E:Rthetauno}, we integrate by parts. For $\sigma<0$
we obtain the following
\[
2e^{\pi i
s/4}\pi^{-s/2}\Gamma(s/2)\Rzeta(s)=\Bigl.\frac{2\tau^{s/2}}{
s}\vartheta(\tau)\Bigr|_{\tau=-1}^{-1+i\infty}-\int_{-1}^{-1+i\infty}
\frac{2\tau^{s/2}}{s}\vartheta_3'(\tau)\,d\tau.
\]
by the known behavior of 
$\vartheta_3(-1+it)$, we get \eqref{E:Rzetatres}. Since the right-hand side defines
an entire function of $s$ the equality is true for all $s$. The poles of $\Gamma(s/2)$, except $s=0$,  
are zeros of $\Rzeta(s)$. 
\end{proof}

\begin{proposition}
For $s=\frac12+it$ with $t$ real, we have
\begin{equation}\label{E:RzetaTheta4}
-s e^{-\pi i s/4}\pi^{-s/2}\Gamma(s/2)\Rzeta(s)=\int_0^\infty g(x,t)e^{if(x,t)}\,dx
=\int_{-1}^{-1+i\infty}(-\tau)^{\frac14+i\frac{t}{2}}
\vartheta_3'(\tau)\,d\tau.
\end{equation}
Here, $f$ and $g$ are real functions defined by
\begin{equation}
g(x,t)=(1+x^2)^{\frac18}e^{\frac{t}{2}\arctan x}\psi(x),\qquad
f(x,t)=\frac{t}{4}\log(1+x^2)-\frac14\arctan x,
\end{equation}
and $\psi(x)=i\vartheta_3'(-1+ix)$ so that
\begin{equation}
\psi(x)=-2\pi\sum_{n=1}^\infty (-1)^n n^2 e^{-\pi n^2x}=\frac{1}{2x^{5/2}}
\sum_{n=0}^\infty((2n+1)^2 \pi-2x)e^{-\frac{\pi}{4x}(2n+1)^2}.
\end{equation}
\end{proposition}

\begin{proof}
By \eqref{E:Rzetatres} we have the following
\[-s\pi^{-s/2}\Gamma(s/2)\Rzeta(s)=e^{-\pi i s/4}
\int_{-1}^{-1+i\infty} \tau^{s/2}\vartheta_3'(\tau)\,d\tau.\]
With  $\tau=-1+iy$ we obtain
\begin{multline*}
\frac{s}{2}\log\tau=\Bigl(\frac{1}{4}+i\frac{t}{2}\Bigr)\Bigl(\frac{1}{2}\log(1+y^2)+
\frac{\pi i}{2}+i\arctan\frac{1}{y}\Bigr)\\
=\frac{1}{8}\log(1+y^2)-\frac{\pi t}{4}-\frac{t}{2}\arctan\frac{1}{y}+
i\Bigl(\frac{t}{4}\log(1+y^2)+\frac{\pi}{8}+\frac{1}{4}\arctan\frac{1}{y}\Bigr)
\end{multline*}
So that 
$-s\pi^{-s/2}\Gamma(s/2)\Rzeta(s)$ is equal to 
\begin{multline*}
e^{-\frac{\pi i}{8} +\frac{\pi t}{4}}\int_0^{+\infty}(1+y^2)^{1/8} e^{-\frac{\pi t}{4}}
e^{-\frac{t}{2}\arctan\frac{1}{y}}\cdot\\
\cdot\exp\Bigl\{i\Bigl(\frac{t}{4}\log(1+y^2)+ \frac{\pi}{8}
+\frac{1}{4}\arctan\frac{1}{y}\Bigr)\Bigr\}\vartheta'(-1+iy) i\,dy
\end{multline*}
Or 
\begin{equation}
\int_0^{+\infty}(1+y^2)^{1/8} e^{-\frac{t}{2}\arctan\frac{1}{y}}
i \vartheta_3'(-1+iy)
\exp\Bigl\{\frac{it}{4}\log(1+y^2)+
\frac{i}{4}\arctan\frac{1}{y}\Bigr\}\,dy.
\end{equation}
Since $\arctan\frac{1}{y}=\frac{\pi}{2}-\arctan y$ we have
\begin{multline*}
-s\pi^{-s/2}\Gamma(s/2)\Rzeta(s)=
\exp\Bigl(\frac{\pi i s}{4}\Bigr)\int_0^{+\infty}(1+y^2)^{1/8} e^{\frac{t}{2}\arctan y}
\psi(y)\cdot\\
\cdot\exp\Bigl\{\frac{it}{4}\log(1+y^2)-
\frac{i}{4}\arctan y\Bigr\}\,dy.
\end{multline*}
This is the first  integral in \eqref{E:RzetaTheta4}; the second integral representation in 
 \eqref{E:RzetaTheta4} is clearly 
equal to the first. 
\end{proof}

\section{A class of functions}\label{S:7}

Lavrik \cite{Lavrik} proved 
\begin{proposition}
For $\tau\in\C$ with $\Re z\ge 0$
\begin{equation}
\Lzeta(\tau,s)=-\frac{\pi^{s/2}\tau^{s/2}}{s\Gamma(s/2)}+
\sum_{n=1}^\infty\frac{\Gamma\bigl(s/2, \pi
n^2\tau\bigr)}{\Gamma(s/2)}\frac{1}{ n^{s}}.
\end{equation}
Then
\begin{equation}
\zeta(s)=\Lzeta(\tau,s)+\chi(s)\Lzeta(1/\tau,1-s).
\end{equation}
\end{proposition}
In particular, \eqref{E:Rzetaserie} implies that  $\Lzeta(i,s)=\Rzeta(s)$.

For $1/\tau=\overline{\tau}$ i.e. when $|\tau|=1$, and taking
$s=1/2+it$ on the critical line, we have
\[\overline{\Lzeta}(\tau,
1-s)=\overline{\Lzeta(\tau,\overline{1-s}) }=\Lzeta(1/\tau,1-s).\]
So  we have 
\begin{equation}
\zeta(s)=\Lzeta(\tau,s)+\chi(s)\overline{\Lzeta}(\tau,1-s),\qquad |\tau|=1,\ \Re\tau\ge0,
\end{equation}
which mentions only a holomorphic function $\Lzeta$ and its conjugate.

We characterize these functions in the next Proposition.

\begin{proposition}\label{P:rzetatype}
The meromorphic functions  $f(s)$ on $\C$ such that
\begin{equation}\label{E:zetaf}
\zeta(s)=f(s)+\chi(s)\overline{f}(1-s)
\end{equation} 
are just those functions that can be written as 
\begin{equation}\label{E:fh}
f(s)=\frac12\zeta(s)+i\frac{\pi^{s/2}}{\Gamma(s/2)}h(s),
\end{equation}
where  $h$ is a meromorphic function on the plane taking real values on the critical line.
\end{proposition}

\begin{proof}
Assume  $f$ can be written as \eqref{E:fh}, then
\begin{multline}\label{E:fhzeta}
\pi^{-s/2}\Gamma(s/2)f(s)+
\pi^{-(1-s)/2}\Gamma((1-s)/2)\overline{f(1-\overline{s})}\\=
\frac{1}{2}\pi^{-s/2}\Gamma(s/2)\zeta(s) +
\frac{1}{2}\pi^{-(1-s)/2}\Gamma((1-s)/2)\overline{\zeta(1-\overline{s})}
+ih(s)-i\overline{h(1-\overline{s})}.\end{multline} Since 
$h$ is real on the real axis for $\sigma=1/2$, i.e. for $s$ in the critical line, we have
\[h(s)=\overline{h(1-\overline{s})},\]
and this and the functional equation of $\zeta(s)$ implies
\[\pi^{-s/2}\Gamma(s/2)f(s)+
\pi^{-(1-s)/2}\Gamma((1-s)/2)\overline{f(1-\overline{s})}=
\pi^{-s/2}\Gamma(s/2)\zeta(s).\] Or equivalently
\[\zeta(s)=f(s)+\chi(s)\overline{f}(1-s),\]
for $s$ in the critical line, but since these are meromorphic functions on all the complex plane. 

Assume now that $f$ satisfies \eqref{E:zetaf} and define
$h$ so that \eqref{E:fh} is satisfied. We have to show that $h$ is real on the critical line. 
Substituting \eqref{E:fh} into \[\pi^{-s/2}\Gamma(s/2)f(s)+
\pi^{-(1-s)/2}\Gamma((1-s)/2)\overline{f(1-\overline{s})}\] we obtain \eqref{E:fhzeta}. Noticing that $f$ satisfies
\eqref{E:zetaf} and the functional equation of  $\zeta$,
\eqref{E:fhzeta} transforms into
\[ih(s)=i\overline{h(1-\overline{s})};\]
which implies that $h$ is real on the real line.
\end{proof}

\section{The functions \texorpdfstring{$\lambda(s)$ and $Y(t)$ }{lambda and Y(t)}} \label{S:8}

By Proposition \ref{P:rzetatype} equation \eqref{E:zetaRzeta} implies that there is a meromorphic
function $\lambda(s)$ such that
\begin{equation}\label{D:lambda}
\Rzeta(s)=\tfrac12\bigl(\zeta(s)+i\lambda(s)\bigr),
\end{equation}
such that $Q(s):=\pi^{-s/2}\Gamma(s)\lambda(s)$ is real on the critical line. This means that 
$\lambda(s)$ satisfies the functional equation
\begin{equation}\label{E:lambdafunct}
Q(s)=\overline{Q}(1-s).
\end{equation}

Since $\Rzeta(s)$ is an entire function, $\lambda(s)$ is a meromorphic function with a unique pole at $s=1$ with residue $1$. The function $\lambda(s)$ have zeros not in the critical line. In the x-ray we may see one, not the first, at a height near 130.

\begin{figure}[H]
\begin{center}
\includegraphics[height=0.6\vsize]{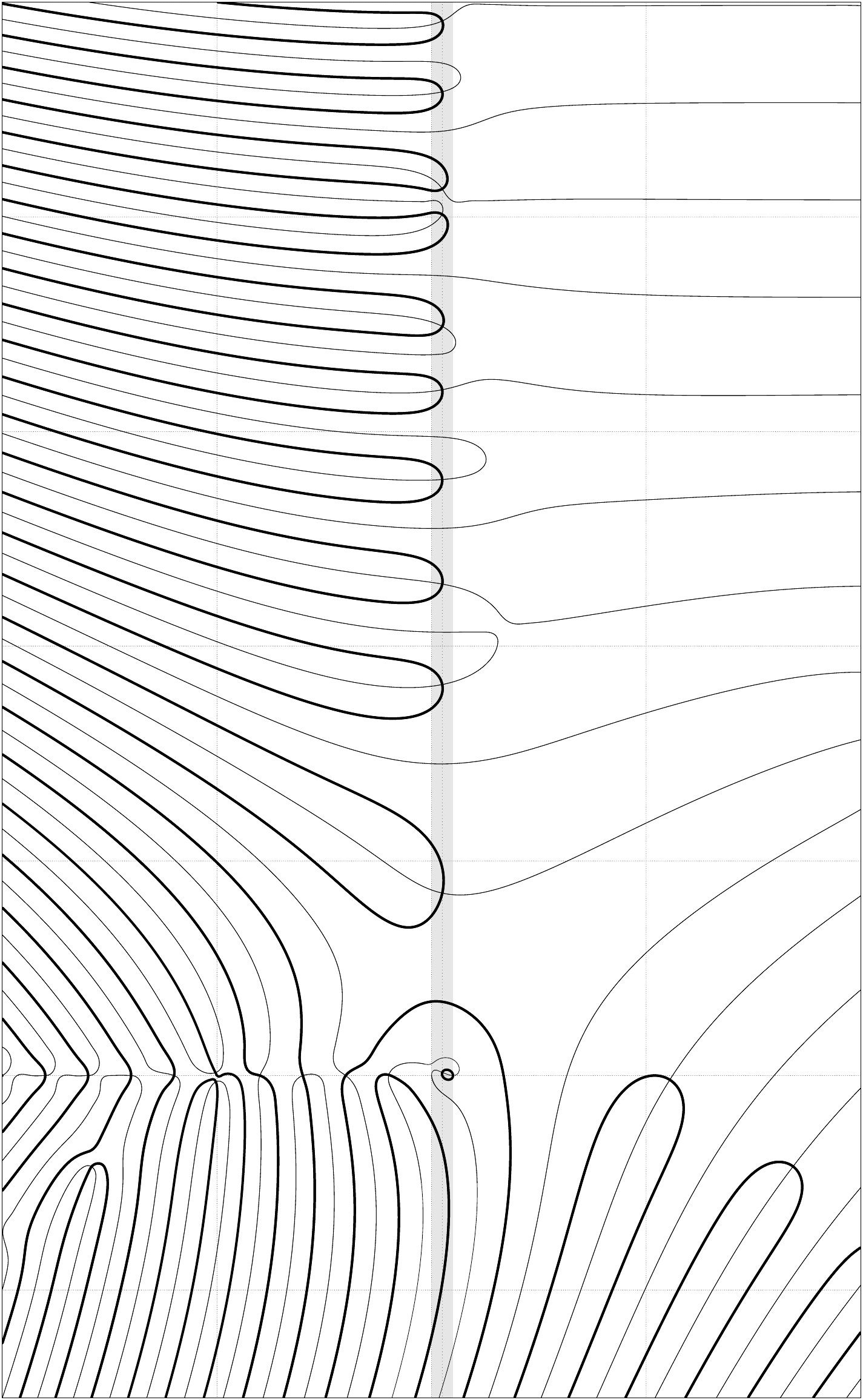}
\includegraphics[height=0.6\vsize]{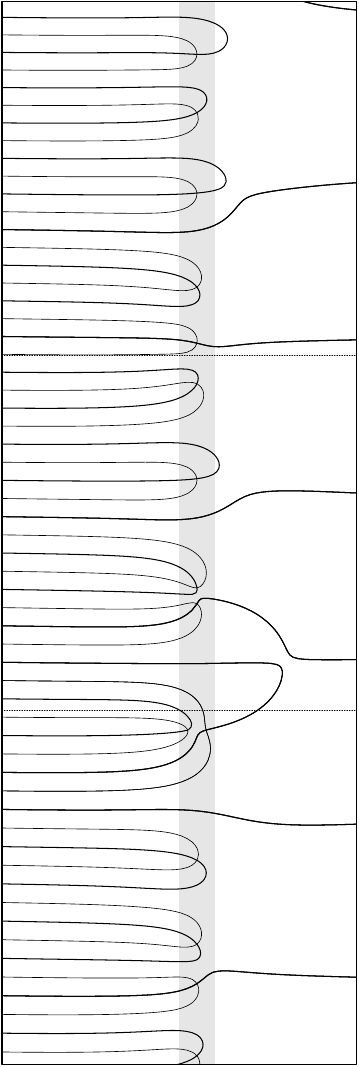}
\caption{Function $\lambda(s)$ for $s\in(-20,20)\times(-15,50)$ and $(-5,5)\times(120,150)$ }
\end{center}
\end{figure}

The function $\vartheta(t)$ is defined as  $\Im\log\Gamma(\frac12+it)-\frac{it}{2}\log \pi$ taking the usual branch of $\log\Gamma(\frac12+it)$. It is a real analytic function. The fact that
$\pi^{-s/2}\Gamma(s/2)\zeta(s)$ is real on the critical line is equivalent to saying that 
$\zeta(\frac12+it)=e^{-i\vartheta(t)}Z(t)$, where $Z(t)$ is a real-valued real analytic function. In the same way, there is a real valued real analytic  function $Y(t)$ such that 
\begin{equation}\label{D:functionY}
\lambda(\tfrac12+it)=e^{-i\vartheta(t)}Y(t).
\end{equation}
It follows that
\begin{equation}\label{E:RzetaZY}
\Rzeta(\tfrac12+it)=\tfrac12 e^{-i\vartheta(t)}\bigl(Z(t)+iY(t)\bigr).
\end{equation}

\end{document}